\documentclass[11pt,a4paper]{amsart}
\usepackage[T1]{fontenc}
\usepackage[utf8]{inputenc}
\usepackage{lmodern,microtype}
\usepackage{amssymb,mathtools}
\usepackage[margin=3cm]{geometry}
\usepackage[hidelinks]{hyperref}

\newtheorem{theorem}{Theorem}[section]
\newtheorem{lemma}[theorem]{Lemma}
\newtheorem{proposition}[theorem]{Proposition}
\newtheorem{mainthm}{Theorem}

\theoremstyle{remark}
\newtheorem{remark}[theorem]{Remark}
\numberwithin{equation}{section}

\newcommand{\R}{\mathbb{R}}
\newcommand{\C}{\mathbb{C}}
\newcommand{\Z}{\mathbb{Z}}
\newcommand{\T}{\mathbb{T}}
\newcommand{\id}{\mathrm{id}}
\DeclareMathOperator{\im}{Im}
\DeclareMathOperator{\re}{Re}
\newcommand{\abs}[1]{\lvert #1\rvert}
\newcommand{\mean}[1]{[#1]}
\newcommand{\ca}{c_{\mathrm A}}
\newcommand{\cR}{c_{\mathrm R}}
\newcommand{\cs}{c_{\ast}}
\newcommand{\Cs}{C_{\ast}}

\title[Linearization of $C^n$ vector fields on the torus]{Linearization of vector fields on the $n$-torus\\ of class $C^n$ with Dini continuous\\ highest derivatives}
\author{Joachim Paulusch}
\thanks{Earlier publications of the author appeared under the name Joachim Albrecht.}

\begin{document}

\begin{abstract}
We consider small perturbations of a constant vector field $\omega$ on the torus $\T^n$, $n\ge2$, where $\omega$ satisfies a Diophantine condition with the smallest possible exponent $n-1$. The perturbation is assumed to be $n$ times continuously differentiable, and the moduli of continuity of its pure derivatives of order $n$ are assumed to satisfy a Dini condition; this is weaker than any H\"older condition. We show that, after a suitable constant vector is added, the perturbed field is $C^1$ conjugate to $\omega$. The proof follows Moser's method of analytic smoothing. In each step it applies an analytic KAM theorem, for which we include a complete proof based on R\"ussmann's optimal estimate for the linear equation. This is a revised English version of a paper published in German in 2005.
\end{abstract}

\maketitle

\section{Introduction}

\subsection{The problem}
Let $\T^n=\R^n/2\pi\Z^n$ and $\omega\in\R^n$. The flow of the constant vector field
\[
  \dot x=\omega
\]
is the simplest flow on the torus: the solutions are $x(t)=x_0+\omega t$, and if the components of $\omega$ are rationally independent, every orbit is dense. We ask whether this picture survives a small perturbation
\[
  \dot x=\omega+f(x),
\]
that is, whether there is a change of coordinates $x=u(\xi)$, with $u-\id$ periodic and small, which turns the perturbed equation back into $\dot\xi=\omega$.

In general this is impossible, even for a constant perturbation $f\equiv c\ne0$: the field $\omega+c$ rotates with a different frequency, and a coordinate change of the above form cannot alter the frequency. Following Moser \cite{Moser66,Moser67}, one therefore admits a \emph{modifying term} $\mu\in\R^n$ and looks for $u$ and $\mu$ such that the coordinate change $x=u(\xi)$ transforms
\[
  \dot x=\omega+f(x)+\mu
\]
into $\dot\xi=\omega$. Since $\dot x=Du(\xi)\dot\xi$, this means
\begin{equation}\label{eq:conj}
  Du(\xi)\,\omega=\omega+f(u(\xi))+\mu\qquad\text{for all }\xi .
\end{equation}
This is a model problem of KAM theory. It shows the analytic difficulties of the Hamiltonian theory of Kolmogorov \cite{Kolmogorov54}, Arnold \cite{Arnold63} and Moser \cite{Moser62}, which was motivated by the stability problem of celestial mechanics, but it avoids most of the geometric work. For introductions to the field we refer to \cite{Arnold89,Moser73,Russmann79}, and to \cite{DiacuHolmes96} for a popular account.

\subsection{Small divisors and the loss of derivatives}
Put $u=\id+v$ in \eqref{eq:conj} and keep only the terms of first order in $v$ and $f$. The result is the linear equation
\begin{equation}\label{eq:linear}
  \partial_\omega v=f+\mu,\qquad \partial_\omega:=\sum_{j=1}^n\omega_j\,\frac{\partial}{\partial\xi_j}.
\end{equation}
If $f=\sum_k\hat f_k\,e^{i\langle k,\xi\rangle}$ is the Fourier series of $f$, then \eqref{eq:linear} is solved by $\mu=-\hat f_0$ and
\[
  v(\xi)=\sum_{k\ne0}\frac{\hat f_k}{i\langle k,\omega\rangle}\,e^{i\langle k,\xi\rangle}.
\]
The \emph{small divisors} $\langle k,\omega\rangle$ come arbitrarily close to zero. To control them, one assumes a \emph{Diophantine condition}
\begin{equation}\label{eq:dioph}
  \abs{\langle k,\omega\rangle}\ge\frac{\gamma}{\abs{k}^{\tau}}\qquad\text{for all }k\in\Z^n\setminus\{0\},
\end{equation}
with constants $\gamma>0$ and $\tau>0$; here $\abs k=\max_j\abs{k_j}$. Such vectors $\omega$ exist only if $\tau\ge n-1$. If $\tau>n-1$, almost every $\omega\in\R^n$ satisfies \eqref{eq:dioph} with some $\gamma>0$. In the limiting case $\tau=n-1$, which is the case of this paper, the set of such $\omega$ has Lebesgue measure zero, but it has the cardinality of the continuum in every open subset of $\R^n$; see \cite{Russmann75} and the references there.

Even under \eqref{eq:dioph}, the solution $v$ is less regular than $f$; roughly speaking, $\tau$ derivatives are lost. For analytic $f$ the loss shows up as a smaller domain of analyticity, and -- in classical KAM-Theory -- the rapid convergence of a Newton iteration has to compensate for it. For finitely differentiable $f$ a Newton iteration would use up the available derivatives after finitely many steps. Moser's remedy \cite{Moser66,Moser67} is to approximate $f$ by analytic functions $f_k$ which are defined on complex strips whose widths $\delta_k$ tend to zero, and to apply the analytic theory to one correction $f_{k+1}-f_k$ at a time. How fast $f_{k+1}-f_k$ tends to zero is governed by the smoothness of $f$.

\subsection{How much smoothness is needed?}
This question has accompanied KAM theory from the beginning. Moser's theorem on invariant curves of area preserving twist maps \cite{Moser62} originally required $333$ derivatives. R\"ussmann reduced this to five derivatives \cite{Russmann70} and then further \cite{Russmann83}. Herman obtained the theorem for perturbations of class $C^p$, $p>3$, and showed by counterexamples that less than three derivatives do not suffice \cite{Herman83,Herman86}.

For a single invariant torus with Diophantine frequency in a Hamiltonian system with $n$ degrees of freedom, Salamon \cite{Salamon04} showed, building on ideas of Moser and P\"oschel \cite{Poschel80}, that a Hamiltonian of class $C^\ell$ with $\ell>2\tau+2$ suffices. For $\tau=n-1$ the H\"older condition beyond $C^{2n}$ can be replaced by a Dini condition on the derivatives of order $2n$, provided the unperturbed system is analytic and integrable \cite{Albrecht07}. Tong and Li \cite{TongLi25} proved results in terms of moduli of continuity for $\tau > n-1$ and for unperturbed systems which need be neither analytic nor integrable. In the other direction, Cheng and Wang \cite{ChengWang13} showed that a Lagrangian torus with a given rotation vector of an integrable system can be destroyed by perturbations which are arbitrarily small in $C^{2n-\varepsilon}$. The persistence of a set of invariant tori of positive measure for finitely differentiable systems was studied by P\"oschel \cite{Poschel82}, Bounemoura \cite{Bounemoura20} and Koudjinan \cite{Koudjinan20}. A comprehensive treatment of the analytic theory of invariant tori in non-degenerate nearly integrable Hamiltonian systems is given by R\"ussmann \cite{Russmann01}.

For vector fields on the torus, the natural threshold is $\tau+1$ derivatives, which equals $n$ in the case $\tau=n-1$. Herman's counterexamples show that perturbations of class $C^{n-\varepsilon}$ do not suffice; see \cite{Herman86} and the discussion in \cite{Poschel21}. P\"oschel \cite{Poschel21} proved that smallness of the perturbation in a certain norm, which measures $\tau+1$ derivatives in an $L^2$ sense, suffices; in particular, the derivatives of order $n$ need not be continuous. His proof uses a slowly convergent iteration in the spirit of R\"ussmann \cite{Russmann10}. The function space used by P\"oschel was studied further by Tong and Li \cite{TongLi26}. For other classes of functions and other arithmetic conditions see, for example, \cite{BounemouraFejoz19,Bounemoura22,TongLi23}.

\subsection{Result}
For a continuous, $2\pi$-periodic function $g\colon\R^n\to\R^n$ let
\[
  K_g(\delta):=\sup_{\abs{x-y}\le\delta}\abs{g(x)-g(y)},\qquad\delta\ge0,
\]
denote its modulus of continuity; here and below $\abs{\cdot}$ is the maximum norm.

\begin{mainthm}\label{thm:main}
Let $n\ge2$, and let $\omega\in\R^n$ satisfy \eqref{eq:dioph} with $\tau=n-1$. Let $f\colon\R^n\to\R^n$ be $2\pi$-periodic in each variable and $n$ times continuously differentiable, and suppose that the moduli of continuity $K_j$ of the derivatives $\partial^nf/\partial x_j^n$ satisfy the Dini condition
\begin{equation}\label{eq:dini}
  \int_0^1\frac{K_j(\delta)}{\delta}\,d\delta<\infty,\qquad j=1,\dots,n.
\end{equation}
Then there is a constant $c_1>0$, depending only on $n$, $\gamma$ and $K_1+\dots+K_n$, with the following property. If $\sup_{\R^n}\abs f\le c_1$, then there are a vector $\mu\in\R^n$ and a $C^1$ diffeomorphism $u$ of $\R^n$, with $u-\id$ $2\pi$-periodic in each variable, such that the coordinate change $x=u(\xi)$ transforms $\dot x=\omega+f(x)+\mu$ into $\dot\xi=\omega$. Moreover, $\abs{Du-E}\le\frac13$ on $\R^n$. In particular, $t\mapsto u(\omega t+\xi_0)$ is a quasi-periodic solution of $\dot x=\omega+f(x)+\mu$ for every $\xi_0\in\R^n$.
\end{mainthm}

Since $u-\id$ is periodic, $u$ induces a $C^1$ diffeomorphism of $\T^n$. Every H\"older condition on the $n$-th derivatives implies \eqref{eq:dini}. The converse is false: a modulus of continuity which behaves like $(\log 1/\delta)^{-a}$ with $a>1$ near zero satisfies \eqref{eq:dini}, but no H\"older condition (cf. \cite{AlbrechtDiss, Albrecht07}). Condition \eqref{eq:dini} is usually called the Dini condition, after Dini's criterion for the convergence of Fourier series.

\subsection{Relation to earlier work}
Theorem~\ref{thm:main} appeared, in German, in \cite{Albrecht05F}. The present paper is a revised version: the proof has been reorganized, the analytic theorem used in it is now proved in full, and the literature has been updated. The existence statement of Theorem~\ref{thm:main} can also be derived from the later and more general theorem of P\"oschel \cite{Poschel21}, because \eqref{eq:dini} implies the smallness condition used there; this is shown in Section~\ref{sec:poschel}. The two proofs are different. The proof given here applies a standard analytic KAM theorem, a proof scheme following Moser \cite{Moser66}, and it is easy to see where the Dini condition enters.

\subsection{Idea of the proof}
Suppose that a transformation $u_k$ solves the problem for the analytic approximation $f_k$. In the coordinates given by $u_k$, the equation with $f_{k+1}$ in place of $f_k$ differs from $\dot\xi=\omega$ only by a term of the size of $f_{k+1}-f_k$, that is, of order $K(\delta_k)\,\delta_k^{\,n}$, where $K$ is a common majorant of the moduli $K_j$. The analytic theorem removes this term by a transformation $v_k$. The derivative of $v_k$ differs from the identity by at most a constant times the size of the term divided by $\delta_k^{\,n}$, that is, by a constant times $K(\delta_k)$: the powers of $\delta_k$ cancel exactly. By the chain rule, the derivative of $u_{k+1}=u_k\circ v_k$ is a product of such factors. It stays under control if $\sum_kK(\delta_k)$ is small, and for a geometric sequence $\delta_k$ this sum is comparable to the integral in \eqref{eq:dini}.

Section~\ref{sec:tools} collects the notation and the auxiliary results. Section~\ref{sec:analytic} contains the analytic theorem and its proof. Theorem~\ref{thm:main} is proved in Section~\ref{sec:proof}, and Section~\ref{sec:poschel} compares the assumptions with those of \cite{Poschel21}.

\section{Notation and auxiliary results}\label{sec:tools}

\subsection{Norms and function spaces}
For $x\in\C^n$ we write $\abs x=\max_j\abs{x_j}$. For a matrix $A=(a_{jl})\in\C^{n\times n}$ we use the corresponding operator norm $\abs A=\max_j\sum_l\abs{a_{jl}}$, so that $\abs{Ax}\le\abs A\abs x$ and $\abs{AB}\le\abs A\abs B$. The identity matrix is $E$. For a bounded function $g$ on a set $M$ we write $\abs g_M=\sup_{x\in M}\abs{g(x)}$. For $s>0$ let
\[
  S(s)=\{x\in\C^n:\ \abs{\im x}<s\}
\]
be the complex strip of width $s$ around $\R^n$, and write $\abs g_s:=\abs g_{S(s)}$. We denote by $P_N(s)$ the space of bounded analytic functions $g\colon S(s)\to\C^N$ which are $2\pi$-periodic in each variable and real for real arguments; matrix valued functions of this kind form the space $P_{n\times n}(s)$. For $g\in P_N(s)$,
\[
  \mean g:=\frac{1}{(2\pi)^n}\int_{[0,2\pi]^n}g(x)\,dx\in\R^N
\]
is the mean value of $g$. Clearly $\abs{\mean g}\le\abs g_s$.

\subsection{Two elementary estimates}
We use the following two facts throughout; their proofs are standard.

\begin{lemma}\label{lem:cauchy}
Let $0<\sigma<s$ and let $F\colon S(s)\to\C^m$ be bounded and analytic.
\begin{enumerate}
\item \emph{(Cauchy's estimate)} $\abs{DF}_{s-\sigma}\le\abs F_s/\sigma$, where $DF$ is the Jacobian matrix.
\item \emph{(Mean value inequality)} If $a,b$ lie in a convex set $M\subseteq S(s)$, then $\abs{F(b)-F(a)}\le\abs{DF}_M\,\abs{b-a}$.
\end{enumerate}
\end{lemma}

\begin{proof}
(1) Let $z\in S(s-\sigma)$ and $h\in\C^n$ with $\abs h\le1$. The function $t\mapsto F(z+th)$ is analytic and bounded by $\abs F_s$ for $\abs t<\sigma$, so Cauchy's inequality in one variable gives $\abs{DF(z)h}\le\abs F_s/\sigma$. Taking the supremum over $h$ gives the operator norm. (2) follows from $F(b)-F(a)=\int_0^1DF(a+t(b-a))\,(b-a)\,dt$.
\end{proof}

Strips are convex, so (2) applies with $M=S(s')$, $s'\le s$. A frequent application is the following. If $w\colon S(s)\to\C^n$ is analytic, real for real arguments, and $\abs{Dw-E}_s\le\frac13$, then for $\xi\in S(s)$ %the point $\re\xi$ is real and
\begin{equation}\label{eq:strip}
  \abs{\im w(\xi)}=\abs{\im\bigl(w(\xi)-w(\re\xi)\bigr)}\le\abs{w(\xi)-w(\re\xi)}\le\tfrac43\,\abs{\im\xi}<\tfrac43\,s,
\end{equation}
because $\abs{Dw}_s\le\frac43$; thus $w$ maps $S(s)$ into $S(\frac43s)$.

\subsection{The common modulus of continuity}
In the situation of Theorem~\ref{thm:main} we put
\begin{equation}\label{eq:K}
  K(\delta):=K_1(\delta)+\dots+K_n(\delta)+\delta,\qquad\delta\ge0.
\end{equation}
The summand $\delta$ only serves to make $K$ positive. The function $K$ majorizes each $K_j$, it is increasing, $K(\delta)>0$ for $\delta>0$, and $K(\delta)\to0$ for $\delta\to0$, since the derivatives $\partial^nf/\partial x_j^n$ are periodic and hence uniformly continuous. Moreover,
\begin{equation}\label{eq:Kdouble}
  K(2\delta)\le2K(\delta),
\end{equation}
because $\abs{g(x)-g(y)}\le\abs{g(x)-g(m)}+\abs{g(m)-g(y)}$ for the midpoint $m$ of $x$ and $y$. By \eqref{eq:dini},
\begin{equation}\label{eq:diniK}
  \int_0^1\frac{K(\delta)}{\delta}\,d\delta<\infty .
\end{equation}

\subsection{Analytic approximation}
The following approximation theorem is proved in \cite[Satz~2.23]{AlbrechtDiss}. It is stated there for real valued functions; for vector valued functions one applies it to each component, which is possible because we use the maximum norm.

\begin{theorem}[Approximation]\label{thm:approx}
Let $m\ge1$, and let $f\colon\R^n\to\R^n$ be $2\pi$-periodic in each variable and $m$ times continuously differentiable. Let $K$ be an increasing function which majorizes the moduli of continuity of $\partial^mf/\partial x_j^m$, $j=1,\dots,n$. Let $(\varrho_k)_{k\ge0}$ be a decreasing sequence of positive numbers with $\varrho_0\le1$ and $\varrho_k\to0$. Then there are a constant $\ca>0$, depending only on $n$ and $m$, and functions $f_k\in P_n(\varrho_k)$ such that
\begin{align*}
  \abs{f-f_k}_{\R^n}&\le\ca\,K(\varrho_k)\,\varrho_k^{\,m}, && k\ge0,\\
  \abs{f_k-f_{k-1}}_{\varrho_k}&\le\ca\,K(\varrho_{k-1})\,\varrho_{k-1}^{\,m}, && k\ge1,\\
  \abs{f_k}_{\varrho_k}&\le\ca\,\abs f_{\R^n}, && k\ge0.
\end{align*}
\end{theorem}

The theorem expresses the classical fact that the smoothness of a function is reflected in how well it can be approximated by analytic functions; compare Jackson's theorem on trigonometric approximation \cite{Achieser67} and Moser's smoothing technique \cite{Moser66,Moser70}.

\subsection{R\"ussmann's estimate for the linear equation}
The analytic theorem rests on the following estimate for the linear equation \eqref{eq:linear}, due to R\"ussmann \cite[Theorem~1.1]{Russmann75}. Elementary estimates of the Fourier series lead to a larger power of $\sigma$ \cite[Lemma~1]{Moser66}. R\"ussmann's power $\sigma^{-\tau}$ is optimal; for related estimates of sums with small divisors see \cite{Russmann76}.

\begin{lemma}[R\"ussmann]\label{lem:russmann}
Let $\omega$ satisfy \eqref{eq:dioph} with $\tau\ge n-1$, $\gamma > 0$, and let $0<\sigma<s\le1$. For every $\varphi\in P_n(s)$ with $\mean\varphi=0$ there is exactly one $w\in P_n(s-\sigma)$ with $\mean w=0$ and $\partial_\omega w=\varphi$, and
\[
  \abs w_{s-\sigma}\le\frac{\cR}{\gamma\,\sigma^{\tau}}\,\abs\varphi_s,
\]
where the constant $\cR>0$ depends only on $n$ and $\tau$. (R\"ussmann measures $k$ with the $\ell^1$-norm; since $\abs k\le\sum_j\abs{k_j}\le n\abs k$, changing the norm in \eqref{eq:dioph} only changes $\gamma$ by a factor $n^{\tau}$, which can be absorbed into $\cR$.)
\end{lemma}

\section{The analytic theorem}\label{sec:analytic}

\subsection{Statement}
The main step in the proof of Theorem~\ref{thm:main} is the following theorem for analytic perturbations. A modifying term $P(x)\lambda$ with a variable matrix $P$ occurs because our iteration works in transformed coordinates, where a constant vector in the original coordinates is no longer constant.

\begin{theorem}[Analytic KAM theorem]\label{thm:analytic}
Let $n\ge2$ and let $\omega$ satisfy \eqref{eq:dioph} with $\tau\ge n-1$. There are constants $\cs>0$ and $\Cs>0$, depending only on $n$, $\gamma$ and $\tau$, with the following property. Let $0<\delta\le1$, and let $g\in P_n(\delta)$ and $P\in P_{n\times n}(\delta)$ with
\begin{equation}\label{eq:Phyp}
  \abs{P-E}_\delta\le\tfrac12 .
\end{equation}
If
\[
  M:=\abs g_\delta\le\cs\,\delta^{\tau+1},
\]
then there are $\lambda\in\R^n$ and $v\colon S(\delta/2)\to S(\delta)$ with $\hat v:=v-\id\in P_n(\delta/2)$ such that
\begin{equation}\label{eq:vconj}
  Dv(\xi)\,\omega=\omega+g(v(\xi))+P(v(\xi))\,\lambda\qquad\text{for all }\xi\in S(\delta/2),
\end{equation}
that is, $x=v(\xi)$ transforms $\dot x=\omega+g(x)+P(x)\lambda$ into $\dot\xi=\omega$. Moreover,
\begin{equation}\label{eq:vest}
  \abs\lambda\le5M\qquad\text{and}\qquad
  \frac{\abs{\hat v}_{\delta/2}}{\delta}+\abs{D\hat v}_{\delta/2}\le\Cs\,\frac{M}{\delta^{\tau+1}} .
\end{equation}
\end{theorem}

\begin{remark}
Theorem~\ref{thm:analytic} is essentially Moser's theorem \cite[pp.~516--517]{Moser66}, together with his extension to a variable coefficient of the parameter \cite[Remark~7, pp.~529--530]{Moser66}. Moser's version has $\delta^{\tau+2}$ instead of $\delta^{\tau+1}$; the improvement comes from Lemma~\ref{lem:russmann}. Since we need the theorem in exactly this form, including the case $\tau=n-1$, we give the proof. It is a Newton iteration of the usual kind; for an alternative approach to this model problem see \cite{Poschel11}.
\end{remark}

The exponent $\tau+1$ in \eqref{eq:vest} has a simple meaning: by Lemma~\ref{lem:russmann}, the solution of the linear equation loses the factor $\sigma^{-\tau}$, and Cauchy's estimate for its derivative costs one more factor $\sigma^{-1}$. In the proof of Theorem~\ref{thm:main} this exponent has to match the order $n=\tau+1$ of the derivatives in Theorem~\ref{thm:approx}; this is the reason why we need the optimal estimate of Lemma~\ref{lem:russmann}.

\subsection{One Newton step}
Let us first explain the idea. Suppose that after some steps we have arrived at the family
\begin{equation}\label{eq:family}
  \dot y=\omega+h(y)+Q(y)\lambda,\qquad\lambda\in\R^n,
\end{equation}
with a small $h$ and a matrix $Q$ close to $E$. We fix part of the parameter, $\lambda=\lambda'+\tilde\lambda$, and look for a transformation $y=\psi(z)=z+w(z)$ which removes $h$ up to a much smaller error. If we neglect all products of small quantities, the transformed equation is $\dot z=\omega+h+Q\lambda'-\partial_\omega w+Q\tilde\lambda$. So we choose $\lambda'$ such that $h+Q\lambda'$ has mean value zero, and we solve the linear equation $\partial_\omega w=h+Q\lambda'$ with Lemma~\ref{lem:russmann}. The error which remains is quadratic in the size of $h$. The following lemma makes this precise.

\begin{lemma}[Newton step]\label{lem:step}
Let $\omega$ satisfy \eqref{eq:dioph} with $\tau\ge n-1$, let $0<s\le1$ and $0<\sigma<s/2$. Let $h\in P_n(s)$ and $Q\in P_{n\times n}(s)$ with $\abs{Q-E}_s\le\frac34$, put $\eta:=\abs h_s$ and
\[
  \varepsilon:=\frac{8\cR}{\gamma}\,\frac{\eta}{\sigma^{\tau+1}},
\]
and assume $\varepsilon\le\frac18$. Then there are $\lambda'\in\R^n$ and $w\in P_n(s-\sigma)$ with the following properties. The map $\psi:=\id+w$ maps $S(s-2\sigma)$ into $S(s-\sigma)$, and with
\[
  Q^+:=(D\psi)^{-1}(Q\circ\psi),\qquad h^+:=(D\psi)^{-1}\bigl((h\circ\psi-h)+(Q\circ\psi-Q)\lambda'\bigr)
\]
we have, on $S(s-2\sigma)$ and for all $\tilde\lambda\in\R^n$,
\begin{equation}\label{eq:stepconj}
  D\psi\,\bigl(\omega+h^++Q^+\tilde\lambda\bigr)=\omega+h\circ\psi+(Q\circ\psi)(\lambda'+\tilde\lambda),
\end{equation}
that is, $y=\psi(z)$ transforms \eqref{eq:family} with $\lambda=\lambda'+\tilde\lambda$ into $\dot z=\omega+h^+(z)+Q^+(z)\tilde\lambda$. Moreover, $h^+\in P_n(s-2\sigma)$, $Q^+\in P_{n\times n}(s-2\sigma)$, and
\begin{gather}
  \abs{\lambda'}\le4\eta,\qquad \abs w_{s-\sigma}\le\varepsilon\sigma,\qquad\abs{Dw}_{s-2\sigma}\le\varepsilon,\label{eq:step1}\\
  \abs{h^+}_{s-2\sigma}\le5\,\varepsilon\,\eta,\qquad \abs{Q^+-Q}_{s-2\sigma}\le3\,\varepsilon .\label{eq:step2}
\end{gather}
\end{lemma}

\begin{proof}
\emph{Choice of $\lambda'$ and $w$.} The mean value $\mean Q$ is a real matrix with $\abs{\mean Q-E}\le\abs{Q-E}_s\le\frac34$. By the Neumann series it is invertible with $\abs{\mean Q^{-1}}\le4$. We put $\lambda':=-\mean Q^{-1}\mean h\in\R^n$; then $\abs{\lambda'}\le4\eta$, and $\varphi:=h+Q\lambda'$ has mean value zero and satisfies $\abs\varphi_s\le\eta+\frac74\cdot4\eta=8\eta$. Lemma~\ref{lem:russmann} gives $w\in P_n(s-\sigma)$ with
\begin{equation}\label{eq:homol}
  \partial_\omega w=Dw\,\omega=h+Q\lambda',\qquad
  \abs w_{s-\sigma}\le\frac{\cR}{\gamma\sigma^\tau}\,8\eta=\varepsilon\sigma .
\end{equation}
Cauchy's estimate gives $\abs{Dw}_{s-2\sigma}\le\varepsilon$. This proves \eqref{eq:step1}.

\emph{The map $\psi$.} For $z\in S(s-2\sigma)$ we have $\abs{\im\psi(z)}\le\abs{\im z}+\abs{w(z)}<s-2\sigma+\sigma/8$, so $\psi$ maps $S(s-2\sigma)$ into $S(s-\sigma)$; the segment from $z$ to $\psi(z)$ lies in $S(s-\sigma)$ as well. Since $\abs{Dw}\le\frac18$, the matrix $D\psi=E+Dw$ is invertible on $S(s-2\sigma)$, with
\begin{equation}\label{eq:Dpsiinv}
  \abs{(D\psi)^{-1}}\le\tfrac87,\qquad\abs{(D\psi)^{-1}-E}\le\frac{\varepsilon}{1-\varepsilon}\le\tfrac87\,\varepsilon .
\end{equation}
Therefore $h^+$ and $Q^+$ are well defined; they are analytic, bounded, periodic and real for real arguments.

\emph{The transformed equation.} By the definition of $Q^+$, $D\psi\,Q^+\tilde\lambda=(Q\circ\psi)\tilde\lambda$. By \eqref{eq:homol}, $D\psi\,\omega=\omega+Dw\,\omega=\omega+h+Q\lambda'$, and by the definition of $h^+$,
\[
  D\psi\,h^+=h\circ\psi+(Q\circ\psi)\lambda'-h-Q\lambda'.
\]
Adding these three identities gives \eqref{eq:stepconj}.

\emph{Estimate of $h^+$.} By Cauchy's estimate $\abs{Dh}_{s-\sigma}\le\eta/\sigma$, and the mean value inequality on $S(s-\sigma)$ gives $\abs{h\circ\psi-h}_{s-2\sigma}\le(\eta/\sigma)\,\varepsilon\sigma=\varepsilon\eta$. Similarly, the function $F:=(Q-E)\lambda'$ satisfies $\abs F_s\le\frac34\abs{\lambda'}\le3\eta$, so that $\abs{(Q\circ\psi-Q)\lambda'}_{s-2\sigma}=\abs{F\circ\psi-F}_{s-2\sigma}\le(3\eta/\sigma)\,\varepsilon\sigma=3\varepsilon\eta$. With \eqref{eq:Dpsiinv} we obtain $\abs{h^+}_{s-2\sigma}\le\frac87\cdot4\varepsilon\eta\le5\varepsilon\eta$.

\emph{Estimate of $Q^+$.} We write $Q^+-Q=\bigl((D\psi)^{-1}-E\bigr)(Q\circ\psi)+(Q\circ\psi-Q)$. The first term is bounded by $\frac87\varepsilon\cdot\frac74=2\varepsilon$. For the second, apply the argument just used to $(Q-E)a$ with $a\in\C^n$, $\abs a\le1$; this gives $\abs{Q\circ\psi-Q}_{s-2\sigma}\le\frac34\varepsilon$. Hence $\abs{Q^+-Q}_{s-2\sigma}\le3\varepsilon$.
\end{proof}

\subsection{Proof of Theorem~\ref{thm:analytic}}
We apply Lemma~\ref{lem:step} repeatedly. The widths of the strips shrink geometrically, and the total loss of width is $\delta/2$:
\[
  \sigma_j:=\frac{\delta}{2^{j+3}},\qquad s_0:=\delta,\qquad s_{j+1}:=s_j-2\sigma_j,\qquad\text{so that } s_j\downarrow\frac\delta2 .
\]
Note that $\sigma_j\le\delta/8<s_j/2$ and $s_j\le1$. We start with $h_0:=g$ and $Q_0:=P$. As long as the assumptions of Lemma~\ref{lem:step} hold, the lemma with $s=s_j$, $\sigma=\sigma_j$, $h=h_j$, $Q=Q_j$ yields $\lambda'_j$, $w_j$ and $\psi_j=\id+w_j$, and we put $h_{j+1}:=h^+$, $Q_{j+1}:=Q^+$. Let
\[
  \eta_j:=\abs{h_j}_{s_j},\qquad \varepsilon_j:=\frac{8\cR}{\gamma}\,\frac{\eta_j}{\sigma_j^{\tau+1}},\qquad A:=5\cdot2^{\tau+1}.
\]

\emph{Quadratic convergence.} By \eqref{eq:step2}, $\eta_{j+1}\le5\varepsilon_j\eta_j$. Since $\sigma_{j+1}=\sigma_j/2$, this gives
\[
  \varepsilon_{j+1}\le\frac{8\cR}{\gamma}\,\frac{5\varepsilon_j\eta_j\,2^{\tau+1}}{\sigma_j^{\tau+1}}=A\,\varepsilon_j^2,\qquad\text{that is,}\qquad A\varepsilon_{j+1}\le(A\varepsilon_j)^2 .
\]
We now choose the constants of Theorem~\ref{thm:analytic}:
\[
  \cs:=\frac{\gamma}{2A\cdot8^{\tau+2}\,\cR},\qquad\Cs:=\frac{4\cdot8^{\tau+2}\,\cR}{\gamma}.
\]
Note that $\cs\Cs=2/A=2^{-\tau}/5$ is a pure number. Then $M\le\cs\delta^{\tau+1}$ means exactly $A\varepsilon_0\le\frac12$, because $\varepsilon_0=8^{\tau+2}\cR M/(\gamma\delta^{\tau+1})$. By induction, $A\varepsilon_{j+1}\le(A\varepsilon_j)^2\le\frac12A\varepsilon_j$, hence
\begin{equation}\label{eq:epsj}
  \varepsilon_{j+1}\le\tfrac12\varepsilon_j,\qquad\sum_{j\ge0}\varepsilon_j\le2\varepsilon_0,\qquad\varepsilon_j\le\varepsilon_0\le\frac1{2A}\le\frac1{40},
\end{equation}
where we used $A\ge20$, which holds because $\tau\ge n-1\ge1$. In particular $\varepsilon_j\le\frac18$ and $\eta_{j+1}\le5\varepsilon_j\eta_j\le\frac18\eta_j$. Moreover, by \eqref{eq:Phyp} and \eqref{eq:step2},
\[
  \abs{Q_j-E}_{s_j}\le\tfrac12+3\sum_{i<j}\varepsilon_i\le\tfrac12+6\varepsilon_0\le\tfrac34 .
\]
Thus the assumptions of Lemma~\ref{lem:step} are satisfied in every step, and the iteration can be continued indefinitely.

\emph{The transformations.} Let $v_0:=\id$ and $v_{j+1}:=v_j\circ\psi_j$; since $\psi_j$ maps $S(s_{j+1})$ into $S(s_j)$, the map $v_j$ is defined on $S(s_j)$, and $v_{j+1}(S(s_{j+1}))\subseteq v_j(S(s_j))\subseteq\dots\subseteq v_0(S(s_0))=S(\delta)$. Hence $g\circ v_j$ and $P\circ v_j$ are defined on $S(s_j)$. Put $\Lambda_j:=\lambda'_0+\dots+\lambda'_{j-1}$. We claim that, on $S(s_j)$,
\begin{equation}\label{eq:vjconj}
  Dv_j\,(\omega+h_j)=\omega+g\circ v_j+(P\circ v_j)\,\Lambda_j .
\end{equation}
For $j=0$ this is trivial. More generally, we show by induction that $Dv_j(\omega+h_j+Q_j\tilde\lambda)=\omega+g\circ v_j+(P\circ v_j)(\Lambda_j+\tilde\lambda)$ for all $\tilde\lambda\in\R^n$. Assume this for $j$, evaluate it at $\psi_j(z)$ with $\lambda'_j+\tilde\lambda$ in place of $\tilde\lambda$, and insert \eqref{eq:stepconj}:
\begin{align*}
  Dv_{j+1}\,(\omega+h_{j+1}+Q_{j+1}\tilde\lambda)
  &=(Dv_j\circ\psi_j)\,D\psi_j\,(\omega+h_{j+1}+Q_{j+1}\tilde\lambda)\\
  &=(Dv_j\circ\psi_j)\bigl(\omega+h_j\circ\psi_j+(Q_j\circ\psi_j)(\lambda'_j+\tilde\lambda)\bigr)\\
  &=\omega+g\circ v_{j+1}+(P\circ v_{j+1})(\Lambda_{j+1}+\tilde\lambda).
\end{align*}
With $\tilde\lambda=0$ this is \eqref{eq:vjconj}.

\emph{Estimates.} We have $\abs{D\psi_j-E}_{s_{j+1}}\le\varepsilon_j$ by \eqref{eq:step1}. From $Dv_{j+1}-E=(Dv_j\circ\psi_j-E)\,D\psi_j+(D\psi_j-E)$ we get, by induction,
\[
  \abs{Dv_j-E}_{s_j}\le\prod_{i<j}(1+\varepsilon_i)-1\le e^{2\varepsilon_0}-1\le3\varepsilon_0,
  \qquad\abs{Dv_j}_{s_j}\le2 ;
\]
here we used $\prod(1+\varepsilon_i)\le\exp\sum\varepsilon_i$, \eqref{eq:epsj}, and $e^x-1\le xe^x$. The segment from $z\in S(s_{j+1})$ to $\psi_j(z)$ lies in $S(s_j-\sigma_j)$, so the mean value inequality gives
\[
  \abs{v_{j+1}-v_j}_{s_{j+1}}=\abs{v_j\circ\psi_j-v_j}_{s_{j+1}}\le\abs{Dv_j}_{s_j}\,\abs{w_j}_{s_j-\sigma_j}\le2\varepsilon_j\sigma_j .
\]
As $\sigma_j\le\sigma_0=\delta/8$, we obtain $\sum_j\abs{v_{j+1}-v_j}_{\delta/2}\le2\sigma_0\sum_j\varepsilon_j\le\frac12\varepsilon_0\delta$. Finally, $\abs{\lambda'_j}\le4\eta_j$ and $\eta_j\le8^{-j}M$, so that $\sum_j\abs{\lambda'_j}\le\frac{32}7M\le5M$.

\emph{Passage to the limit.} The maps $v_j$ are all defined on $S(\delta/2)$ and converge uniformly there to a map $v$ with $\abs{v-\id}_{\delta/2}\le\frac12\varepsilon_0\delta$. Since the $v_j$ are analytic, so is $v$, and $Dv_j\to Dv$ pointwise on $S(\delta/2)$ by Weierstrass' theorem; hence $\abs{Dv-E}_{\delta/2}\le3\varepsilon_0$. The map $\hat v=v-\id$ is periodic and real for real arguments, because all $v_j-\id$ are. For $\xi\in S(\delta/2)$ we have $\abs{\im v(\xi)}<\frac\delta2+\frac12\varepsilon_0\delta<\delta$, and likewise for $v_j(\xi)$. The sequence $\Lambda_j$ converges to some $\lambda\in\R^n$ with $\abs\lambda\le5M$, and $\abs{Dv_j\,h_j}_{\delta/2}\le2\eta_j\to0$. Letting $j\to\infty$ in \eqref{eq:vjconj} gives \eqref{eq:vconj}. Finally,
\[
  \frac{\abs{\hat v}_{\delta/2}}{\delta}+\abs{D\hat v}_{\delta/2}\le\tfrac12\varepsilon_0+3\varepsilon_0\le4\varepsilon_0=\Cs\,\frac{M}{\delta^{\tau+1}} .
\]
This completes the proof of Theorem~\ref{thm:analytic}.\qed

\section{Proof of Theorem~\ref{thm:main}}\label{sec:proof}

From now on $\tau=n-1$, so that $\tau+1=n$, and $K$ is defined by \eqref{eq:K}.

\subsection{The scales}
Let $\delta_k:=2^{-k}\delta_0$, where $\delta_0\in(0,\frac12]$ will be chosen in Lemma~\ref{lem:theta}. We apply Theorem~\ref{thm:approx} with $m=n$ and $\varrho_k:=2\delta_k$; note that $\varrho_0\le1$ and $S(\varrho_{k+1})=S(\delta_k)$. Using \eqref{eq:Kdouble}, we obtain functions $f_k\in P_n(2\delta_k)$ with
\begin{align}
  \abs{f-f_k}_{\R^n}&\le2^{n+1}\ca\,K(\delta_k)\,\delta_k^{\,n},\label{eq:ap1}\\
  \abs{f_{k+1}-f_k}_{\delta_k}&\le2^{n+1}\ca\,K(\delta_k)\,\delta_k^{\,n}
  \le4^{n+1}\ca\,K(\delta_{k+1})\,\delta_{k+1}^{\,n},\label{eq:ap2}\\
  \abs{f_0}_{\delta_0}&\le\ca\,\abs f_{\R^n}.\label{eq:ap3}
\end{align}
In the second inequality of \eqref{eq:ap2} we used $\delta_k=2\delta_{k+1}$ and \eqref{eq:Kdouble} once more.

\subsection{Where the Dini condition enters}
Before we set up the iteration, we explain which quantity has to be controlled. Suppose that a transformation $u_k$ solves the problem for $f_k$, and that $\abs{(Du_k)^{-1}}\le\frac32$. To pass to $f_{k+1}$, we shall apply Theorem~\ref{thm:analytic} on the strip $S(\delta_{k+1})$ to the perturbation
\begin{equation}\label{eq:g}
  g_{k+1}:=(Du_k)^{-1}\,(f_{k+1}-f_k)\circ u_k .
\end{equation}
By \eqref{eq:ap2} its size is at most $\frac32\cdot4^{n+1}\ca\,K(\delta_{k+1})\,\delta_{k+1}^{\,n}$. Theorem~\ref{thm:analytic} requires that this be at most $\cs\delta_{k+1}^{\,n}$, so we measure it in units of $\cs\delta_{k+1}^{\,n}$ and define
\begin{equation}\label{eq:theta}
  \theta_k:=\frac{3\cdot4^{n+1}\ca}{2\,\cs}\,K(\delta_k),\qquad k\ge0 .
\end{equation}
Then $\abs{g_{k+1}}_{\delta_{k+1}}\le\cs\,\theta_{k+1}\,\delta_{k+1}^{\,n}$, and Theorem~\ref{thm:analytic} applies if $\theta_{k+1}\le1$. By \eqref{eq:vest}, the resulting transformation $v$ satisfies $\abs{Dv-E}_{\delta_{k+2}}\le\Cs\cs\,\theta_{k+1}$. Here the powers of $\delta_{k+1}$ have cancelled. The transformations are composed, so their derivatives are multiplied, and
\[
  \prod_k\bigl(1+\Cs\cs\,\theta_k\bigr)\le\exp\Bigl(\Cs\cs\sum_k\theta_k\Bigr).
\]
Hence the derivatives remain under control if the sum of the $\theta_k$ is small. The following lemma shows that this is the case for small $\delta_0$, and this is the only place where the Dini condition is used.

\begin{lemma}\label{lem:theta}
The sequence $(\theta_k)$ is decreasing, and $\delta_0\in(0,\frac12]$ can be chosen, depending only on $n$, $\gamma$ and $K$, such that
\begin{equation}\label{eq:thetasum}
  \theta_0\le1\qquad\text{and}\qquad\Cs\cs\sum_{k=0}^\infty\theta_k\le\log\frac43 .
\end{equation}
\end{lemma}

\begin{proof}
The sequence decreases because $K$ is increasing. For the same reason, for $k\ge1$,
\[
  K(\delta_k)\log2=K(\delta_k)\int_{\delta_k}^{\delta_{k-1}}\frac{ds}{s}\le\int_{\delta_k}^{\delta_{k-1}}\frac{K(s)}{s}\,ds,
\]
and summation over $k$ gives
\[
  \sum_{k=0}^\infty K(\delta_k)\le K(\delta_0)+\frac{1}{\log2}\int_0^{\delta_0}\frac{K(s)}{s}\,ds .
\]
By \eqref{eq:diniK} and since $K(\delta)\to0$, the right hand side tends to zero as $\delta_0\to0$. The constants in \eqref{eq:theta} and \eqref{eq:thetasum} depend only on $n$ and $\gamma$.
\end{proof}

\subsection{The iteration}
We fix $\delta_0$ as in Lemma~\ref{lem:theta} and put
\[
  c_1:=K(\delta_0)\,\delta_0^{\,n}>0 .
\]

\begin{proposition}\label{prop:iter}
Let $\abs f_{\R^n}\le c_1$. Then for every $k\ge0$ there are $\mu_k\in\R^n$ and a map $u_k\colon S(\delta_{k+1})\to S(\delta_k)$ with $u_k-\id\in P_n(\delta_{k+1})$ such that, on $S(\delta_{k+1})$,
\begin{align}
  Du_k\,\omega&=\omega+f_k\circ u_k+\mu_k,\label{eq:i1}\\
  \abs{Du_k-E}_{\delta_{k+1}}&\le\exp\Bigl(\Cs\cs\sum_{l=0}^k\theta_l\Bigr)-1\le\frac13 .\label{eq:i2}
\end{align}
Moreover, with $u_{-1}:=\id$ and $\mu_{-1}:=0$,
\begin{align}
  \abs{u_k-u_{k-1}}_{\delta_{k+1}}&\le\tfrac43\,\Cs\cs\,\theta_k\,\delta_k,\label{eq:i3}\\
  \abs{\mu_k-\mu_{k-1}}&\le5\,\cs\,\theta_k\,\delta_k^{\,n}.\label{eq:i4}
\end{align}
\end{proposition}

Equation \eqref{eq:i1} says that $x=u_k(\xi)$ transforms $\dot x=\omega+f_k(x)+\mu_k$ into $\dot\xi=\omega$.

\begin{proof}
We use induction on $k$. Once $u_k-\id\in P_n(\delta_{k+1})$ and \eqref{eq:i2} are known, \eqref{eq:strip} shows that $u_k$ maps $S(\delta_{k+1})$ into $S(\frac43\delta_{k+1})\subseteq S(\delta_k)$. So this property need not be checked separately.

\emph{The case $k=0$.} By \eqref{eq:ap3}, the choice of $c_1$, and \eqref{eq:theta},
\[
  \abs{f_0}_{\delta_0}\le\ca\,c_1=\ca\,K(\delta_0)\,\delta_0^{\,n}\le\cs\,\theta_0\,\delta_0^{\,n}.
\]
Since $\theta_0\le1$, Theorem~\ref{thm:analytic} applies with $g=f_0$, $P=E$ and $\delta=\delta_0$. It gives $\lambda$ and $v$, and we put $u_0:=v$, $\mu_0:=\lambda$. Then \eqref{eq:i1} is \eqref{eq:vconj}, and \eqref{eq:vest} gives $\abs{\mu_0}\le5\cs\theta_0\delta_0^n$, $\abs{u_0-\id}_{\delta_1}\le\Cs\cs\theta_0\delta_0$ and $\abs{Du_0-E}_{\delta_1}\le\Cs\cs\theta_0\le e^{\Cs\cs\theta_0}-1$. The bound $\frac13$ in \eqref{eq:i2} follows from \eqref{eq:thetasum}, because $e^{\log(4/3)}-1=\frac13$.

\emph{From $k$ to $k+1$.} Let $U:=(Du_k)^{-1}$ on $S(\delta_{k+1})$. By \eqref{eq:i2} and the Neumann series,
\[
  \abs U_{\delta_{k+1}}\le\tfrac32,\qquad\abs{U-E}_{\delta_{k+1}}=\abs{U(E-Du_k)}_{\delta_{k+1}}\le\tfrac32\cdot\tfrac13=\tfrac12 .
\]
Since $u_k$ maps $S(\delta_{k+1})$ into $S(\delta_k)$, where $f_{k+1}$ and $f_k$ are defined, the function $g_{k+1}=U\,(f_{k+1}-f_k)\circ u_k$ of \eqref{eq:g} belongs to $P_n(\delta_{k+1})$; periodicity holds because $u_k-\id$ is periodic. As shown in the preceding subsection, $\abs{g_{k+1}}_{\delta_{k+1}}\le\cs\theta_{k+1}\delta_{k+1}^n$ with $\theta_{k+1}\le\theta_0\le1$. Hence Theorem~\ref{thm:analytic} applies with $g=g_{k+1}$, $P=U$ and $\delta=\delta_{k+1}$. It gives $\lambda\in\R^n$ and $v\colon S(\delta_{k+2})\to S(\delta_{k+1})$ with
\begin{equation}\label{eq:veq}
  Dv\,\omega=\omega+g_{k+1}\circ v+(U\circ v)\,\lambda\qquad\text{on }S(\delta_{k+2})
\end{equation}
and, by \eqref{eq:vest},
\begin{equation}\label{eq:vk}
  \abs\lambda\le5\cs\theta_{k+1}\delta_{k+1}^{\,n},\qquad
  \abs{v-\id}_{\delta_{k+2}}\le\Cs\cs\theta_{k+1}\delta_{k+1},\qquad
  \abs{Dv-E}_{\delta_{k+2}}\le\Cs\cs\theta_{k+1}.
\end{equation}
We put $u_{k+1}:=u_k\circ v$ and $\mu_{k+1}:=\mu_k+\lambda$. Then $u_{k+1}-\id$ is analytic, bounded, periodic and real for real arguments.

To verify \eqref{eq:i1}, we multiply \eqref{eq:veq} by $Du_k\circ v$. On the left we obtain $Du_{k+1}\,\omega$. On the right, $(Du_k\circ v)(U\circ v)=E$, and by \eqref{eq:i1} at the point $v(\xi)\in S(\delta_{k+1})$,
\begin{align*}
  Du_{k+1}\,\omega&=(Du_k\circ v)\,\omega+(f_{k+1}-f_k)\circ u_{k+1}+\lambda\\
  &=\omega+f_k\circ u_{k+1}+\mu_k+(f_{k+1}-f_k)\circ u_{k+1}+\lambda
  =\omega+f_{k+1}\circ u_{k+1}+\mu_{k+1}.
\end{align*}
Estimate \eqref{eq:i4} follows from \eqref{eq:vk}. For \eqref{eq:i2} we write
\[
  Du_{k+1}-E=(Du_k\circ v-E)\,Dv+(Dv-E).
\]
With $a:=\exp\bigl(\Cs\cs\sum_{l\le k}\theta_l\bigr)-1$ and $b:=\Cs\cs\theta_{k+1}$ this gives $\abs{Du_{k+1}-E}_{\delta_{k+2}}\le a(1+b)+b=(1+a)(1+b)-1$, which is at most $\exp\bigl(\Cs\cs\sum_{l\le k+1}\theta_l\bigr)-1\le\frac13$. For \eqref{eq:i3} we note that the segment from $\xi\in S(\delta_{k+2})$ to $v(\xi)\in S(\delta_{k+1})$ lies in the convex set $S(\delta_{k+1})$, on which $\abs{Du_k}\le\frac43$. The mean value inequality gives
\[
  \abs{u_{k+1}(\xi)-u_k(\xi)}=\abs{u_k(v(\xi))-u_k(\xi)}\le\tfrac43\,\abs{v(\xi)-\xi}\le\tfrac43\,\Cs\cs\,\theta_{k+1}\,\delta_{k+1}.
\]
This completes the induction.
\end{proof}

\subsection{Convergence}
\begin{proof}[Proof of Theorem~\ref{thm:main}]
By \eqref{eq:i4} and Lemma~\ref{lem:theta}, $\sum_k\abs{\mu_k-\mu_{k-1}}<\infty$, so $\mu_k\to\mu\in\R^n$. By \eqref{eq:i3}, the maps $u_k$ converge uniformly on $\R^n$ to a continuous map $u$, and $u-\id$ is periodic.

The derivatives converge as well. Let $x\in\R^n$. The ball $\{z\in\C^n:\abs{z-x}<\delta_{k+1}\}$ lies in $S(\delta_{k+1})$, so Cauchy's estimate, applied as in the proof of Lemma~\ref{lem:cauchy} with radius $\delta_{k+1}$, and \eqref{eq:i3} give
\[
  \abs{Du_k(x)-Du_{k-1}(x)}\le\frac{\abs{u_k-u_{k-1}}_{\delta_{k+1}}}{\delta_{k+1}}\le\tfrac83\,\Cs\cs\,\theta_k,\qquad k\ge0 .
\]
The right hand side is summable by Lemma~\ref{lem:theta}. Hence $(Du_k)$ converges uniformly on $\R^n$, $u$ is continuously differentiable, $Du=\lim Du_k$, and $\abs{Du-E}_{\R^n}\le\frac13$ by \eqref{eq:i2}.

Consequently $u-\id$ is periodic and Lipschitz continuous with constant at most $\frac13$. Such a map is a diffeomorphism: it is injective, since $\abs{u(x)-u(y)}\ge\frac23\abs{x-y}$; it is surjective,  
since for given $y$ the equation $u(x)=y$ is equivalent to the fixed point equation $x=y-(u(x)-x)$, whose right hand side is a contraction of $\mathbb R^n$ with constant $\frac13$; by Banach's fixed point theorem it has a (unique) solution.
The inverse of $u-\id$ is $C^1$ by the inverse function theorem, since $Du$ is invertible.

Finally, $f_k\to f$ uniformly on $\R^n$ by \eqref{eq:ap1}, and $f$ is uniformly continuous. Hence $f_k\circ u_k\to f\circ u$ uniformly on $\R^n$, and letting $k\to\infty$ in \eqref{eq:i1} at real points gives
\[
  Du\,\omega=\omega+f\circ u+\mu \qquad\text{on }\R^n,
\]
which is \eqref{eq:conj}. Thus $x=u(\xi)$ transforms $\dot x=\omega+f(x)+\mu$ into $\dot\xi=\omega$, and $t\mapsto u(\omega t+\xi_0)$ is a solution of $\dot x=\omega+f(x)+\mu$, since its derivative is $Du(\omega t+\xi_0)\,\omega$.
\end{proof}

\begin{remark}
The same proof works if $\omega$ satisfies \eqref{eq:dioph} with an integer $\tau\ge n-1$, if $f$ is of class $C^{\tau+1}$, and if the pure derivatives of order $\tau+1$ satisfy the Dini condition \eqref{eq:dini}. One applies Theorem~\ref{thm:approx} with $m=\tau+1$ and replaces $n$ by $\tau+1$ in all exponents of $\delta_k$; the powers of $\delta_k$ again cancel. Taking $\tau=n$, one obtains a result for almost all frequency vectors $\omega$, at the price of one more derivative.
\end{remark}

\section{Comparison with P\"oschel's smallness condition}\label{sec:poschel}

Write $f=\sum_k\hat f_k\,e^{i\langle k,x\rangle}$ and fix an integer $b\ge2$. For $\tau=n-1$, the theorem of P\"oschel \cite{Poschel21} applies if the quantity
\[
  \|f\|_{n,b}:=\sum_{\nu\ge0}\Bigl(\sum_{b^{\nu-1}<\|k\|\le b^\nu}\abs{\hat f_k}_2^2\,\|k\|^{2n}\Bigr)^{1/2}
\]
is sufficiently small; here $\abs{\cdot}_2$ is the Euclidean norm, the sum excludes $k=0$ (the mean value of $f$ can be absorbed into the modifying term), and $\|k\|$ is a norm on $\Z^n$. Which norm is used only affects constants; we assume $\abs k\le\|k\|\le n\abs k$, which covers the $\ell^1$-, $\ell^2$- and maximum norms. We show that the assumptions of Theorem~\ref{thm:main} imply that $\|f\|_{n,b}$ is small if $\abs f_{\R^n}$ is small. P\"oschel's theorem then provides a constant modifying term and a $C^1$ transformation conjugating the perturbed field to $\omega$. Hence the existence statement of Theorem~\ref{thm:main} can also be derived from \cite{Poschel21}; as in Theorem~\ref{thm:main}, the admissible size of $\abs f_{\R^n}$ then depends on the moduli $K_j$. A related observation, formulated with an $L^2$ modulus of continuity, is made in \cite[Theorem~1.6]{TongLi26}. Throughout this section, $C$ denotes constants which depend only on $n$ and $b$.

\emph{Large frequencies.} Let $g_j:=\partial^nf/\partial x_j^n$; its Fourier coefficients are $(ik_j)^n\hat f_k$. Fix $\nu$ with $b^{\nu-1}\ge4n$ and put $N:=\lfloor b^{\nu-1}/(2n)\rfloor\ge1$. For every $k$ with $\|k\|>b^{\nu-1}$ there is an index $j$ with $\abs{k_j}=\abs k\ge\|k\|/n>2N$, and then $\|k\|^{2n}\le n^{2n}\abs{k_j}^{2n}$. Hence
\begin{equation}\label{eq:block}
  \sum_{b^{\nu-1}<\|k\|\le b^\nu}\abs{\hat f_k}_2^2\,\|k\|^{2n}\le n^{2n}\sum_{j=1}^n\ \sum_{\abs{k_j}\ge2N}\abs{\hat f_k}_2^2\,\abs{k_j}^{2n}.
\end{equation}
Let $V_N$ be the de la Vall\'ee Poussin mean of order $N$ with respect to the variable $x_j$. It multiplies the Fourier coefficient with index $k$ by a factor which vanishes for $\abs{k_j}\ge2N$. By Parseval's identity, the inner sum on the right of \eqref{eq:block} is at most $(2\pi)^{-n}\int_{[0,2\pi]^n}\abs{g_j-V_Ng_j}_2^2\,dx$. For fixed values of the other variables, $g_j$ is a periodic function of $x_j$ with modulus of continuity at most $K_j$. By Jackson's theorem and the uniform bound for the norms of the de la Vall\'ee Poussin operators \cite{Achieser67}, $\abs{g_j-V_Ng_j}_2\le C\,K_j(1/N)$ pointwise. Since $1/N\le4n\,b^{1-\nu}$, the subadditivity of $K_j$ gives $K_j(1/N)\le C\,K_j(b^{-\nu})$. Together, the $\nu$-th term of $\|f\|_{n,b}$ is bounded by $C\sum_jK_j(b^{-\nu})$. As in the proof of Lemma~\ref{lem:theta}, $K_j(b^{-\nu})\log b\le\int_{b^{-\nu}}^{b^{1-\nu}}K_j(s)\,s^{-1}\,ds$, so these terms are summable by \eqref{eq:dini}. Hence the sum over $\nu\ge\nu_0$ is smaller than any given number if $\nu_0$ is large enough, depending only on $n$, $b$ and $K_1,\dots,K_n$.

\emph{Small frequencies.} For each $\nu<\nu_0$, the $\nu$-th term is at most $b^{n\nu_0}\bigl(\sum_k\abs{\hat f_k}_2^2\bigr)^{1/2}\le\sqrt n\,b^{n\nu_0}\abs f_{\R^n}$, again by Parseval's identity. These finitely many terms are small if $\abs f_{\R^n}$ is small.

\subsection*{Acknowledgements}
This revision was prepared with substantial assistance from the AI system Claude (Anthropic; models Claude Sonnet 5.5 and Claude Opus 5.5). It was used to translate and restructure the text, to search the literature, to draft the proof of Theorem~\ref{thm:analytic}, and to check the estimates. The author has verified all arguments and references and takes full responsibility for the content. Any remaining errors remain with the author.


\begin{thebibliography}{99}\small
\bibitem{Achieser67} N.\,I. Achieser, \emph{Vorlesungen \"uber Approximationstheorie}, 2nd ed., Akademie-Verlag, Berlin, 1967.
\bibitem{AlbrechtDiss} J. Albrecht, \emph{\"Uber die Stabilit\"at invarianter Tori in Hamiltonschen Systemen, die bis auf eine endlich oft differenzierbare St\"orung analytisch und integrabel sind}, Dissertation, Johannes Gutenberg-Universit\"at Mainz, 2005.
\bibitem{Albrecht05F} J. Albrecht, Kleine Nenner bei $2\pi$-periodischen Vektorfeldern, in: H.\,P. Delfosse, H.\,R. Yousefi (eds.), \emph{``Wer ist weise? der gute Lehr von jedem annimmt''. Festschrift f\"ur Michael Albrecht zum 65. Geburtstag}, Traugott Bautz, Nordhausen, 2005, 371--388.
\bibitem{Albrecht07} J. Albrecht, On the existence of invariant tori in nearly-integrable Hamiltonian systems with finitely differentiable perturbations, \emph{Regul. Chaotic Dyn.} \textbf{12} (2007), 281--320.
\bibitem{Arnold63} V.\,I. Arnold, Proof of a theorem of A.\,N. Kolmogorov on the invariance of quasi-periodic motions under small perturbations of the Hamiltonian, \emph{Russian Math. Surveys} \textbf{18} (1963), no.~5, 9--36.
\bibitem{Arnold89} V.\,I. Arnold, \emph{Mathematical methods of classical mechanics}, 2nd ed., Springer, New York, 1989.
\bibitem{Bounemoura20} A. Bounemoura, Positive measure of KAM tori for finitely differentiable Hamiltonians, \emph{J. \'Ec. polytech. Math.} \textbf{7} (2020), 1113--1132.
\bibitem{Bounemoura22} A. Bounemoura, Optimal linearization of vector fields on the torus in non-analytic Gevrey classes, \emph{Ann. Inst. H. Poincar\'e C Anal. Non Lin\'eaire} \textbf{39} (2022), 501--528.
\bibitem{BounemouraFejoz19} A. Bounemoura, J. F\'ejoz, KAM, $\alpha$-Gevrey regularity and the $\alpha$-Bruno--R\"ussmann condition, \emph{Ann. Sc. Norm. Super. Pisa Cl. Sci. (5)} \textbf{19} (2019), 1225--1279.
\bibitem{ChengWang13} C.-Q. Cheng, L. Wang, Destruction of Lagrangian torus for positive definite Hamiltonian systems, \emph{Geom. Funct. Anal.} \textbf{23} (2013), 848--866.
\bibitem{DiacuHolmes96} F. Diacu, P. Holmes, \emph{Celestial encounters}, Princeton University Press, Princeton, 1996.
\bibitem{Herman83} M.\,R. Herman, Sur les courbes invariantes par les diff\'eomorphismes de l'anneau, Vol.~1, \emph{Ast\'erisque} \textbf{103--104} (1983).
\bibitem{Herman86} M.\,R. Herman, Sur les courbes invariantes par les diff\'eomorphismes de l'anneau, Vol.~2, \emph{Ast\'erisque} \textbf{144} (1986).
\bibitem{Kolmogorov54} A.\,N. Kolmogorov, On conservation of conditionally periodic motions for a small change in Hamilton's function, \emph{Dokl. Akad. Nauk SSSR} \textbf{98} (1954), 527--530.
\bibitem{Koudjinan20} C.\,E. Koudjinan, A KAM theorem for finitely differentiable Hamiltonian systems, \emph{J. Differential Equations} \textbf{269} (2020), 4720--4750.
\bibitem{Moser62} J. Moser, On invariant curves of area-preserving mappings of an annulus, \emph{Nachr. Akad. Wiss. G\"ottingen Math.-Phys. Kl. II} (1962), 1--20.
\bibitem{Moser66} J. Moser, A rapidly convergent iteration method and non-linear partial differential equations I, II, \emph{Ann. Scuola Norm. Sup. Pisa (3)} \textbf{20} (1966), 265--315, 499--535.
\bibitem{Moser67} J. Moser, Convergent series expansions for quasi-periodic motions, \emph{Math. Ann.} \textbf{169} (1967), 136--176.
\bibitem{Moser70} J. Moser, On the construction of almost periodic solutions for ordinary differential equations, in: \emph{Proc. Internat. Conf. on Functional Analysis and Related Topics (Tokyo, 1969)}, Univ. of Tokyo Press, Tokyo, 1970, 60--67.
\bibitem{Moser73} J. Moser, \emph{Stable and random motions in dynamical systems}, Princeton University Press, Princeton, 1973.
\bibitem{Poschel80} J. P\"oschel, \"Uber invariante Tori in differenzierbaren Hamiltonschen Systemen, \emph{Bonner Math. Schriften} \textbf{120} (1980).
\bibitem{Poschel82} J. P\"oschel, Integrability of Hamiltonian systems on Cantor sets, \emph{Comm. Pure Appl. Math.} \textbf{35} (1982), 653--696.
\bibitem{Poschel11} J. P\"oschel, KAM \`a la R, \emph{Regul. Chaotic Dyn.} \textbf{16} (2011), 17--23.
\bibitem{Poschel21} J. P\"oschel, KAM below $C^n$, preprint, 2021, arXiv:2104.01866.
\bibitem{Russmann70} H. R\"ussmann, Kleine Nenner I: \"Uber invariante Kurven differenzierbarer Abbildungen eines Kreisringes, \emph{Nachr. Akad. Wiss. G\"ottingen Math.-Phys. Kl. II} (1970), 67--105.
\bibitem{Russmann75} H. R\"ussmann, On optimal estimates for the solutions of linear partial differential equations of first order with constant coefficients on the torus, in: J. Moser (ed.), \emph{Dynamical systems, theory and applications}, Lecture Notes in Physics \textbf{38}, Springer, Berlin, 1975, 598--624.
\bibitem{Russmann76} H. R\"ussmann, Note on sums containing small divisors, \emph{Comm. Pure Appl. Math.} \textbf{29} (1976), 755--758.
\bibitem{Russmann79} H. R\"ussmann, Konvergente Reihenentwicklungen in der St\"orungstheorie der Himmelsmechanik, in: K. Jacobs (ed.), \emph{Selecta Mathematica V}, Heidelberger Taschenb\"ucher \textbf{201}, Springer, Berlin, 1979, 93--260.
\bibitem{Russmann83} H. R\"ussmann, On the existence of invariant curves of twist mappings of an annulus, in: \emph{Geometric dynamics (Rio de Janeiro, 1981)}, Lecture Notes in Math. \textbf{1007}, Springer, Berlin, 1983, 677--718.
\bibitem{Russmann01} H. R\"ussmann, Invariant tori in non-degenerate nearly integrable Hamiltonian systems, \emph{Regul. Chaotic Dyn.} \textbf{6} (2001), 119--204.
\bibitem{Russmann10} H. R\"ussmann, KAM-iteration with nearly infinitely small steps in dynamical systems of polynomial character, \emph{Discrete Contin. Dyn. Syst. Ser. S} \textbf{3} (2010), 683--718.
\bibitem{Salamon04} D.\,A. Salamon, The Kolmogorov--Arnold--Moser theorem, \emph{Math. Phys. Electron. J.} \textbf{10} (2004), Paper 3.
\bibitem{TongLi23} Z. Tong, Y. Li, Towards sharp regularity: full-dimensional tori in $C^\infty$ vector fields over $\T^\infty$, preprint, 2023, arXiv:2306.08211.
\bibitem{TongLi25} Z. Tong, Y. Li, Towards continuity: universal frequency-preserving KAM persistence and remaining regularity in nearly non-integrable Hamiltonian systems, \emph{Commun. Contemp. Math.} \textbf{27} (2025), doi:10.1142/S021919972450038X.
\bibitem{TongLi26} Z. Tong, Y. Li, Sharp regularity of a weighted Sobolev space over $\T^n$ and its relation to finitely differentiable KAM theory, preprint, 2026, arXiv:2604.04665.
\end{thebibliography}
\end{document}